\documentclass[a4paper,12pt]{amsart}
\usepackage{ifthen}
\usepackage[paper=a4paper]{geometry}
\usepackage[english]{babel}
\usepackage[latin1]{inputenc}
\usepackage[T1]{fontenc}
\usepackage{graphicx}
\usepackage{tabls}
\usepackage{multicol}
\usepackage{setspace}
\usepackage{amssymb}	
\usepackage{enumerate}
\usepackage{amsmath}
\usepackage{amsthm}
\usepackage{comment}
\usepackage{mathrsfs}
\usepackage{float}
\usepackage{amsfonts}
\numberwithin{equation}{section}
\nonstopmode
\newcommand{\B}{{\mathcal B}}
\newcommand{\Bp}{{\mathcal B_p}}

\newcommand{\Co}{{\mathcal{C}o}}
\newcommand{\C}{{\mathbb C}}
\newcommand{\D}{{\mathbb D}}

\newcommand{\F}{{\mathcal F}}

\newcommand{\R}{{\mathbb R}}

\newcommand{\es}{{\mathcal S}}

\newcommand{\X}{{\mathbf X}}

\newcommand{\bD}{{\overline{\mathbb D}}}
\newcommand{\sphere}{{\widehat{\mathbb C}}}

\newcommand{\inv}{^{-1}}

\newcommand{\aand}{{\quad\text{and}\quad}}
\newcommand{\phii}{{\varphi}}

\begin{document}
\bibliographystyle{amsplain}
 \theoremstyle{plain}
 \newtheorem{Definition}{Definition}
 \newtheorem{Theorem}[Definition]{Theorem}
 \newtheorem{Proposition}[Definition]{Proposition}
 \newtheorem{Corollary}[Definition]{Corollary}
 \newtheorem{Lemma}[Definition]{Lemma}
 \theoremstyle{definition}
 \newtheorem{Remark}[Definition]{Remark}
 \newtheorem*{acknowledgement}{Acknowledgment}

\newcounter{alphabet}
\newcounter{tmp}
\newenvironment{Thm}[1][]{\refstepcounter{alphabet}%
\bigskip%
\noindent%
{\bf Theorem \Alph{alphabet}}%
\ifthenelse{\equal{#1}{}}{}{ (#1)}%
{\bf .}
\itshape}{\vskip 8pt}
\newcommand{\Ref}[1]{\setcounter{tmp}{\ref{#1}}\Alph{tmp}}

\newenvironment{pf}[1][]{%
 \vskip 3mm
 \noindent
 \ifthenelse{\equal{#1}{}}%
  {{\slshape Proof. }}%
  {{\slshape #1.} }%
 }%
{\qed\bigskip}

%Title etc
\title[Coefficient estimates of analytic endomorphisms]{Coefficient estimates of analytic endomorphisms of
  the unit disk fixing a point with applications to concave functions}
\author[Rintaro Ohno]{Rintaro Ohno}
\address{Graduate School of Information Sciences, \endgraf Tohoku University, Sendai, 980-8579, Japan.}
\email{rohno@ims.is.tohoku.ac.jp}
\thanks{The first author was supported by Grant-in-Aid for JSPS Fellows No. $26 \cdot 2855$.}
\author[Toshiyuki Sugawa]{Toshiyuki Sugawa}
\address{Graduate School of Information Sciences, \endgraf Tohoku University, Sendai, 980-8579, Japan.}
\email{sugawa@math.is.tohoku.ac.jp}
%\date{\today}

%Abstract
\subjclass[2010]{to be specified}
\begin{abstract}
In this note, we discuss the coefficient regions of analytic self-maps
of the unit disk with a prescribed fixed point.
As an application, we solve the Fekete-Szeg\H{o} problem for
normalized concave functions with a prescribed pole in the unit disk.
\end{abstract}

%\dedicatory{Dedicated to the 90th birthday of Professor Yukio Kusunoki}

\maketitle
\thispagestyle{empty}

%Introduction
\section{Introduction}
\label{sec::intro}
Let $\D= \{z \in \C : |z|<1\}$ denote the unit disk in the complex plane $\C.$ 
The class $\B_p$ for $p\in\D$ will mean the set of holomorphic maps
$\phii: \D \to \D$ satisfying $\phii(p)=p.$
In what follows, without loss of generality, we will always assume that
$0\le p<1.$

A function $\phii\in\Bp$ can be expanded about the origin in the form
\begin{equation}\label{eq::phi}
\phii(z)= c_0+c_1z+c_2z^2+\cdots
=\sum_{n=0}^{\infty} c_n z^n.
\end{equation}
Note that $|c_n| \leq 1$ for each $n.$
We define the coefficient body $\X_n(\F)$ of order $n\geq 0$ for a class $\F$ 
of analytic functions at the origin as the set
$$
\big\{(c_0,c_1,\dots,c_n)\in\C^{n+1}: 
\varphi(z)=c_0+c_1z+\dots +c_nz^n+O(z^{n+1})
\text{ for some } \varphi\in\F\big\}.
$$
Note that $\pi_{m,n}(\X_n(\F))=\X_m(\F)$ for $0\le m<n,$
where $\pi_{m,n}:\C^{n+1}\to\C^{m+1}$ is the projection
$(c_0,c_1,\dots,c_n)\mapsto (c_0,c_1,\dots, c_m).$

Obviously, $\X_0(\B_0)=\{0\}$ and $\X_1(\B_0)=\{(0,c): |c|\le1\}.$
In the present paper, we describe $\X_n(\B_p)$ for $n=0,1$ and
$0<p<1.$
Note that the authors describe $X_2(\B_p)$ in \cite{OS15} to
investigate the second Hankel determinant.
In the following, it is convenient to put
$$
P=p+\frac1p=\frac{1+p^2}{p}.
$$
Note that $P>2.$

\begin{Theorem}
\label{thm::c01}
Let $p \in (0,1)$.
\begin{enumerate}
\item[$(\mathrm{i})$] $\X_0 (\B_p) = \left\{c_0\in\C : \left| c_0 - P\inv \right| 
\leq P\inv \right\}.$
For a function $\varphi(z)=c_0+c_1z+\cdots$ in $\B_p,$ $c_0\in\partial\X_0(\B_p)$ 
if and only if $\varphi$ is an analytic automorphism of $\D.$
\medskip
\item[$(\mathrm{ii})$] $\X_1 (\B_p)=\left\{(c_0,c_1)\in\C^2: 
\left| c_1 -(1-Pc_0+c_0^2) \right| 
\le P\left[P^{-2}-\left| c_0 - P\inv \right|^2\right]\right\}.
$
In other words, a pair $(c_0,c_1)$ of complex numbers is contained in 
$\X_1(\B_p)$ if and only if
\begin{equation}\label{eq::c}
c_0=P\inv(1-\sigma_0)
\aand
c_1=P^{-2}\big[1+(P^2-2)\sigma_0+\sigma_0^2\big]+P\inv(1-|\sigma_0|^2)\sigma_1
\end{equation}
for some $\sigma_0, \sigma_1\in\bD.$

Moreover, for a function $\varphi(z)=c_0+c_1z+\cdots$ in $\B_p,$ 
$(c_0,c_1)\in\partial\X_1(\B_p)$ 
if and only if $\varphi$ is either an analytic automorphism of $\D$
or a Blaschke product of degree $2.$
\end{enumerate}
\end{Theorem}

Our motivation of the present study comes from an intimate relation
between $\B_p$ and the class $\Co_p$ of concave functions $f$ normalized by $f(0)=f'(0)-1=0$ with a pole at $p.$
%We will apply this result to a Fekete-Szeg\"o-type problem of concave functions
%in the last section. 
Here, a meromorphic function $f$ on $\D$ is said to be {\it concave}, 
if it maps $\D$ conformally onto a concave domain in the Riemann sphere
$\sphere=\C\cup\{\infty\};$ in other words, $f$ is a univalent meromorphic
function on $\D$ such that $\C\setminus f(\D)$ is convex.
The class $\Co_p$ is intensively studied in recent years by Avkhadiev, Bhowmik,
Pommerenke, Wirths and others 
(see e.g. \cite{APW04, APW06, AW02, AW05, AW07, BPW11}).

The following representation of concave functions in terms of functions
in $\B_p$ is contained in the first author's paper \cite{Ohno13}.

\begin{Thm}
\label{thm::intrep}
Let $0<p<1$ and put $P=p+1/p.$
A meromorphic function $f$ on $\D$
is contained in the class $\mathcal{C}o_p$ if and only if 
there exists a function $\phii\in\B_p$ such that
\begin{equation}
\label{equ::intp}
f'(z)= (1-Pz+z^2)^{-2}
\exp \int_0^z \frac{-2 \phii(\zeta)}{1-\zeta \phii(\zeta)} d\zeta.
\end{equation}
\end{Thm}

For a given function $f\in\Co_p$ with the expansion
\begin{equation}\label{eq::f}
f(z)= z+a_2z^2+a_3z^3+\dots
=\sum_{n=1}^{\infty} a_n z^n,\quad |z|<p,
\end{equation}
we consider the Fekete-Szeg\H{o} functional
$$
\Lambda_\mu(f)=a_3-\mu a_2^2
$$
for a real number $\mu.$
For example, $\Lambda_1(f)=a_3-a_2^2=S_f(0)/6,$ where 
$S_f=(f''/f')'-(f''/f')^2/2$ is the Schwarzian derivative of $f.$
For some background of the Fekete-Szeg\H{o} functional,
the reader may refer to \cite{CKS07} and references therein.
As an application of Theorem \ref{thm::c01}, we will prove the following.

\begin{Theorem}\label{thm::main}
Let $0<p<1$ and $\mu\in\R$ and put $P=p+1/p.$
Then the maximum $\Phi(\mu)$ of the Fekete-Szeg\H{o} functional
$|\Lambda_\mu(f)|$ over $f\in\Co_p$ is given as follows:
$$
\Phi(\mu)=
\begin{cases}
(1-\mu)P^2-1 & \text{if}\quad \mu\le \mu_1(P), \\
%& \text{if}\quad 0\le\mu\le\mu_1, \\
-\dfrac{1}{3} (P^3-2P+3)+\dfrac{(P+2)^2(2P-1)^2}{12(P+3\mu)}
& \text{if}\quad \mu_1(P)\le\mu\le\mu_2(P), \\
%& \text{if}\quad \mu_2\le\mu\le\mu_3, \\
\Psi(P,\mu)
& \text{if}\quad \mu_2(P)\le\mu\le\mu_4(P), \\
(\mu-1)P^2+1 & \text{if}\quad \mu_4(P)\le\mu.
\end{cases}
$$
Here,
$$
\Psi(P,\mu)=\begin{cases}
P^2-3-\mu(P^2-4+4P^{-2}) & \null \\
\qquad\qquad\text{if either}~ P_2\le P\le P_*,~ \mu_3^-(P)\le\mu\le\mu_3^+(P)
&\null \\
\qquad\qquad\text{or}~ P_*\le P,~ \mu_2(P)\le \mu\le \mu_3^+(P), &\null \\
(1-\mu)P(P^2-2)\sqrt{\dfrac{P^2-4\mu}{4\mu\{(1-\mu)(P^2-1)^2-1\}}}
&\text{otherwise,}
\end{cases}
$$
and
\begin{align*}
\mu_1(P)&=\frac12-\frac1{3P}, \\
\mu_2(P)&=\begin{cases}
\dfrac{1}{72}\left(4+P^2+4P^4-\sqrt{16P^8+8P^6-543P^4+1160P^2+16}\right)&\quad \text{if}~ P\le P_*, \\
\dfrac{P(3P+2)}{6(P^2-2)}
&\quad \text{if}~ P_*\le P,
\end{cases} \\
\mu_3^\pm(P)&=\frac{P^2(3P^4-12P^2+14)\pm P^2 \sqrt{P^8-16P^6+84P^4-176P^2+132}}%
{4(P^2-1)(P^2-2)^2} \\
\mu_4(P)&=\frac{3P^4-4P^2-2+\sqrt{P^8-12P^4+16P^2+4}}{4P^2(P^2-1)},
\end{align*}
where $P_*\approx 2.88965$ is the unique root of the polynomial
$$
U(P)=6P^4-P^3-38P^2-28P+4
$$
on the interval $2<P<+\infty$ and $P_2\approx 2.82343$ is the largest root of
the polynomial
$$
V(P)=P^8-16P^6+84P^4-176P^2+132
$$
on the positive real axis.
Moreover, 
$$
\frac13<\mu_1<\frac12<\mu_2<\mu_4<\frac89,
$$
on the interval $2<P,$ and $\mu_2<\mu_3^-<\mu_3^+<\mu_4$ on $P_2<P<P_*$
whereas $\mu_3^-<\mu_2<\mu_3^+<\mu_4$ on $P_*<P.$
\end{Theorem}

We see numerically that $p_*\in(0,1)$ satisfying $P_*=p_*+1/p_*$
is approximately $0.401984.$
Also, we have $p_2\approx 0.415252$ for $p_2\in(0,1)$ with
$P_2=p_2+1/p_2.$

\begin{figure}
\begin{center}
\includegraphics{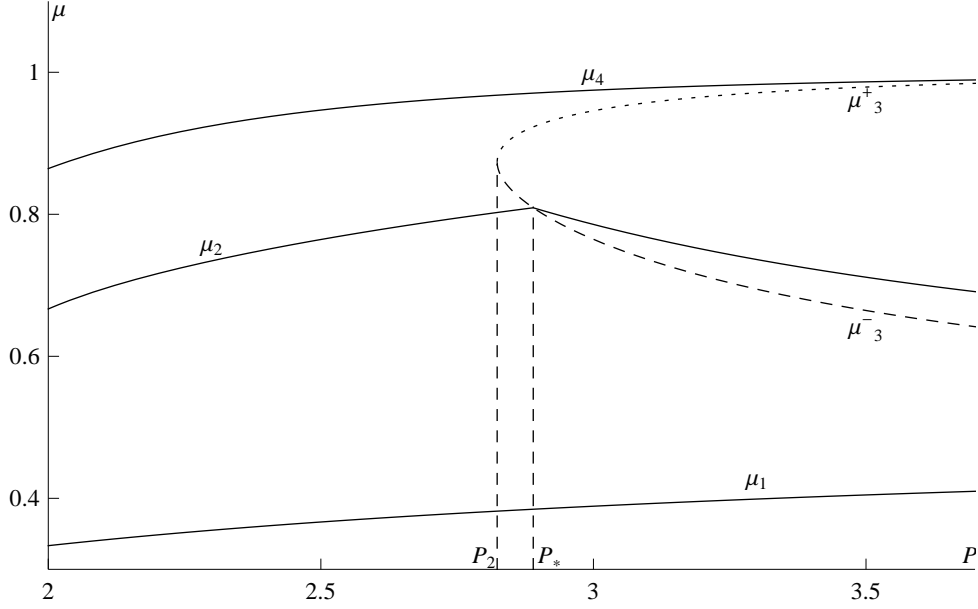}
\caption{The graphs of $\mu_1(P),\mu_2(P),\mu_3^{\pm}(P)$ and $\mu_4(P)$ 
in the $P\mu$-plane.}
\label{fig::Pmu}
\end{center}
\end{figure}

The behaviour of $\mu_1(P),\mu_2(P),\mu_3^{\pm}(P)$ and $\mu_4(P)$ can be observed in Figure \ref{fig::Pmu}.

\medskip
The Fekete-Szeg\H{o} problem was solved by Bhowmik, Ponnusamy and Wirths in 
\cite{BPW11} for the different but related classes $\Co(\alpha)$ for $1<\alpha\le2.$
Here, by definition, $f\in\Co(\alpha)$ if $f\in\es$ with $f(1)=\infty,$
if $\C\setminus f(\D)$ is convex, and if
the opening angle of the image $f(\D)$ at $\infty$ is $\le \pi\alpha.$
It is interesting to observe that
the case $\alpha=2$ of their main theorem in \cite{BPW11} agrees with
the limiting case of our Theorem \ref{thm::main} as $p\to1^-$
(equivalently, $P\to2^+$).

\medskip
With the special choice $\mu=0$, we have the following known fact.

\begin{Corollary}
Let $f(z)=z+a_1z+a_2z^2+\cdots$ be a function in $\mathcal{C}o_p.$
Then the following sharp inequality holds:
$$ |a_3| \leq P^2-1=p^2+1+ \frac{1}{p^2}.$$
\end{Corollary}

Indeed, the above inequality is still valid as long as $f$ is univalent meromorphic
on $\D$ with a pole at $p$ (see Jenkins \cite{Jen62}).
Avkhadiev, Pommerenke and Wirths \cite{APW04} (see also \cite{AW07})
proved the even stronger result that the variability region
of $a_3$ over $f\in\Co_p$ is given as $|a_3-P^2+2|\le 1.$
(This can also be proved by our method given below.)

Since $\Phi(1)=1$ by Theorem \ref{thm::main}, we get another corollary.

\begin{Corollary}
Let $0<p<1$ and suppose that
$f(z)=z+a_1z+a_2z^2+\cdots$ is a function in $\mathcal{C}o_p.$
Then the following sharp inequality holds:
$$ |a_3-a_2^2| \leq 1.$$
\end{Corollary}

Recall that $6(a_3-a_2^2)=S_f(0)$ is the Schwarzian derivative of $f$ evaluated
at $z=0.$ 
The inequality $|a_3-a_2^2|\le 1$ is valid
for a univalent holomorphic function $f(z)=z+a_2z^2+a_3z^3+\cdots$
on $\D$ (see, for instance, \cite[Ex.~1 in p.~70]{Duren:univ}).
Indeed, it is obtained by a simple application of Gronwall's area theorem
for the function $1/f(1/w).$
Since the Schwarzian derivative $S_f$ is unchanged under the post-composition
with M\"obius transformations, the above corollary can also be obtained 
from this classical result.

In the final section, we will focus on the variability region of 
$\Lambda_1(f)=a_3-a_2^2$ over $f\in\Co_p.$
Section 2 will be devoted to the proof of Theorem \ref{thm::c01}.
In order to apply Theorem \ref{thm::c01} to concave functions, in Section 3,
we consider a maximum value problem for a quadratic polynomial over
the closed unit disk.
The proof of Theorem \ref{thm::main} will be given in Section 4.

\section{Proof of Theorem \ref{thm::c01}}
\label{sec::proof}

For the proof of Theorem \ref{thm::c01}, we recall a 
useful lemma due to Dieudonn\'e 
(see \cite[p.~198]{Duren:univ} for instance).
To clarify the equality case in the lemma, we will give a proof briefly.
Throughout this section, it is helpful to use the special automorphism
\begin{equation}\label{eq::Ta}
T_a(z)=\frac{a-z}{1-\bar az}
\end{equation}
of $\D$ for $a\in\D.$
This is indeed an analytic involution of $\D$ and interchanges $0$ and $a.$
Moreover,
$$
T_a'(z)=\frac{|a|^2-1}{(1-\bar az)^2}.
$$
In particular,
$$
T_a'(0)=|a|^2-1
\aand
T_a'(a)=\frac1{|a|^2-1}.
$$

\begin{Lemma}[Dieudonn\'e's Lemma]
\label{lem::dieudonne}
Let $z_0, w_0 \in \D$ with $|w_0|<|z_0|.$
Then the region of values of $w=g'(z_0)$ for
holomorphic functions $g:\D\to\D$ with $g(0)=0$ and $g(z_0)=w_0$
is given as the closed disk
\begin{equation}
\label{equ::dieudonne}
\left|w - \frac{w_0}{z_0} \right| 
\leq \frac{|z_0|^2 - |w_0|^2}{|z_0| (1-|z_0|^2)}.
\end{equation}
Equality holds if and only if
$g$ is a Blaschke product of degree $2$ fixing $0.$
\end{Lemma}

\begin{pf}
The function $h(z)=g(z)/z$ is an analytic endomorphism of $\D$ which
sends $z_0$ to $\omega_0.$
Thus $H=T_{\omega_0}\circ h\circ T_{z_0}$ belongs to $\B_0.$
The Schwarz lemma now gives $|H'(0)|\le 1$ which turns out to be equivalent
to \eqref{equ::dieudonne} with $w=g'(z_0).$
Moreover, equality holds if and only if $H(z)=\zeta z$ for some
$\zeta\in\partial\D,$ which means $h$ is an analytic automorphism of $\D.$
\end{pf}

In view of the proof, we have a concrete form of $g$ in the equality case:
$$
g(z)= z T_{\omega_0}(\zeta T_{z_0}(z))
$$
for a constant $\zeta\in\partial\D,$ where $\omega_0=w_0/z_0\in\D.$

We are now ready to prove Theorem \ref{thm::c01}.

\begin{pf}[Proof of Theorem \ref{thm::c01}]
For a function $\varphi\in\B_p,$ we consider $\psi=T_p\circ\varphi\circ T_p:
\D\to\D.$
Then $\psi\in\B_0.$
The Schwarz lemma implies $|\psi(p)|\le p.$
Namely,
$$
|T_p(c_0)|=\left| \frac{p-c_0}{1-p \, c_0} \right| \leq p,
$$
which is equivalent to
\begin{equation}\label{eq::pc}
0\le |1-pc_0|^2-\left|1 - \dfrac{c_0}p\right|^2
=\dfrac{1-p^4}{p^2}\left[
\left(\dfrac p{1+p^2}\right)^2-\left| c_0 - \dfrac{p}{1+p^2} \right|^2
\right].
\end{equation}
The range is optimal because the function $\varphi$ corresponding
to $\psi(z)=T_p(c_0)z/p$ belongs to $\B_p.$
Suppose now that $c_0\in\partial\X_0(\B_p).$
Then, by the above argument, we have $\psi(z)=\zeta z,$
where $\zeta=T_p(c_0)/p\in\partial\D.$
Thus $\varphi(z)=T_p(\zeta T_p(z))$ is an analytic automorphism
of $\D$ fixing $p.$
Hence the first assertion follows.

For the second assertion, we use Dieudonn\'e's lemma.
Note that
$$
\psi'(p)=T_p'(c_0)\cdot\varphi'(0)\cdot T_p'(p)
=\frac{c_1}{(1-pc_0)^2}.
$$
Applying Dieudonn\'e's lemma to the function $\psi$
with the choices $z_0=p$ and $w_0=\psi(p)=T_p(c_0),$ we get
$$
\left|  \frac{c_1}{(1-p\, c_0)^2} -\frac{p - c_0}{p\, (1-p\, c_0)}\right| 
\leq
\frac{p^2 - \left|\frac{p - c_0}{1-p\, c_0}\right|^2}{p \,(1-p^2)}.
$$
Here, if $|w_0|=p=|z_0|,$ the above inequality (in fact, equality)
holds obviously.
Note that the above range of $c_1$ for a fixed $c_0$ is optimal 
by Dieudonn\'e's lemma.
Using the identity in \eqref{eq::pc}, we obtain the first description of
the set $\X_1(\B_p).$
The second description of $\X_1(\B_p)$
is obtained by letting $\sigma_0=P(P\inv-c_0)=1-Pc_0$
and $\sigma_1=(c_1-(1-Pc_0+c_0^2))/(P\inv-P|c_0-P\inv|^2)
=P(c_1-P^{-2}(1+(P^2-2)\sigma_0+\sigma_0^2))/(1-|\sigma_0|^2).$

We now prove the final assertion.
Suppose that $(c_0,c_1)\in\partial\X_1(\B_p)$
for a function $\varphi(z)=c_0+c_1z+\cdots$ in $\B_p.$
By part (i), we know that $c_0\in\partial\X_0(\B_p)$ if and only if
$\varphi$ is an analytic automorphism of $\D$ fixing $p.$
%Note that $\pi_{0,1}\inv(\partial\X_0(\B_p))\subset \partial\X_1(\B_p).$
%Thus, it is enough to show that $(c_0,c_1)\in\partial\X_1(\B_p)
%\setminus\pi_{0,1}\inv(\partial\X_0(\B_p))$ if and only if
%$\varphi$ is a Blaschke product of degree 2 fixing $p.$
%We first show the ``only if" part.
%We can now assume that $c_0$ and $c_1$ have the representations in
%\eqref{eq::c} with $|w_0|<1$ and $|w_1|=1.$
Thus we may assume that $c_0$ is an interior point of $\X_0(\B_p);$
namely, $|T_p(c_0)|<p.$
Then, by the equality case in Dieudonn\'e's lemma,
$\psi=T_p\circ\varphi\circ T_p$ is a Blaschke product of degree 2 fixing $0.$
Therefore, we have proved the ``only if" part.
The ``if" part is easy to check.
\end{pf}

\section{Maximum value problem for a quadratic polynomial}
In order to apply Theorem \ref{thm::c01} for concave functions, 
we consider the following problem:
What is the value of the quantity
\begin{equation}
\label{equ::max}
Y(a,b,c)=\max_{z\in\bD}
\left(  \left|a + b z +c z^2\right| + 1- |z|^2 \right)
\end{equation}
for real numbers $a,b,c$?

In fact, a more general and symmetric problem was considered in \cite{CKS07}.
Let
$$
\Omega(A,B,K,L,M)=
\max_{u,v\in\bD}\left\{
|A|(1-|u|^2)+|B|(1-|v|^2)+|Ku^2+2Muv+Lv^2|
\right\}
$$
for $A, B, K, L, M\in\C.$
When $K, L, M$ are all real numbers, the value of $\Omega(A,B,K,L,M)$
is computed in \cite[Theorem 3.1]{CKS07}.
By virtue of the maximum modulus principle, one can see that
$$
\Omega(1,0,c,a,b/2)
=\max_{u\in\bD, v\in\partial\D}\left\{
(1-|u|^2)+|cu^2+buv+av^2|\right\}
=Y(a,b,c).
$$
As an immediate consequence of Theorem 3.1 in \cite{CKS07},
we obtain the following result.
(Note that, under the notation adopted in \cite{CKS07},
$\max\{\Phi_1,\Phi_2\}\ge0$ because of $B=0$ so that
$S\ge |A|+|B|=1$ in the case (3c) of Theorem 3.1 in \cite{CKS07}.)

\begin{Proposition}
\label{prop::max}
Let  $Y(a,b,c)$ be the quantity defined in \eqref{equ::max}
for real numbers $a,b,c.$
When $ac\ge0,$
$$
Y(a,b,c)=\begin{cases}
|a|+|b|+|c| &\quad \text{if}~ |b|\ge 2(1-|c|),\\
1+|a|+\dfrac{b^2}{4(1-|c|)} &\quad\text{if}~ |b|< 2(1-|c|).
\end{cases}
$$
When $ac<0,$
\begin{equation}\label{eq::Q}
Y(a,b,c)=\begin{cases}
1-|a|+\dfrac{b^2}{4(1-|c|)} 
&\quad\text{if}~ -4ac(c^{-2}-1)\le b^2 \text{ and } |b|<2(1-|c|), \\
1+|a|+\dfrac{b^2}{4(1+|c|)} 
&\quad\text{if}~ b^2<\min\{4(1+|c|)^2,-4ac(c^{-2}-1)\}, \\
R(a,b,c)
&\quad \text{otherwise,}
\end{cases}
\end{equation}
where
\begin{equation}\label{eq::R}
R(a,b,c)=\begin{cases}
|a|+|b|-|c|&\quad\text{if}~ |c|(|b|+4|a|)\le|ab|, \\
-|a|+|b|+|c|&\quad\text{if}~ |ab|\le |c|(|b|-4|a|), \\
(|c|+|a|)\sqrt{1-\dfrac{b^2}{4ac}} &\quad\text{otherwise}.
\end{cases}
\end{equation}
\end{Proposition}

\section{Proof of Theorem \ref{thm::main}}

Let $p\in (0,1)$ and put $P=p+1/p$ as before.
For a given function $f\in\Co_p$ with expansion \eqref{eq::f},
there is a unique function $\varphi\in\B_p$ with expansion \eqref{eq::phi}
such that the representation formula \eqref{equ::intp} holds.
A straightforward computation yields
$$
a_2= P -c_0
\aand
a_3=P^2- \frac{1}{3}\left( c_1 -c_0^2 +4P c_0 +2 \right).
$$
For $\mu \in \R$, by substituting the expressions in \eqref{eq::c},
we obtain
\begin{align}\label{eq::FS}
&\quad~ a_3-\mu a_2^2
=\frac13\left[(1-3\mu)c_0^2+2(3\mu-2)Pc_0-c_1+(3-\mu)P^2-2\right] \\
&=P^2-2-\mu(P-P\inv)^2+\left(1-2\mu(1-P^{-2})\right)\sigma_0-\mu P^{-2}\sigma_0^2
-\frac{(1-|\sigma_0|^2)\sigma_1}{3P}. \notag
\end{align}
Since $\sigma_1$ is an arbitrary point in $\bD,$ we get the sharp
inequality
$$
|a_3-\mu a_2^2|\le
\frac{1}{3P}\left\{|a+b\sigma_0+c\sigma_0^2|+1-|\sigma_0|^2\right\},
$$
where
$$
a=3P\big[P^2-2-\mu (P-P\inv)^2\big], \quad
b=3P-6\mu(P-P^{-1})
\aand
c=-3\mu P\inv.
$$
Therefore, in terms of the quantity introduced in the last section,
we can express $\Phi(\mu)$ by
$$
\Phi(\mu)=\sup_{f\in\Co_p}\Lambda_\mu(f)
=\frac{1}{3P}Y(a,b,c).
$$

Observe that $a$ changes its sign at $\mu=\mu_a:=(P^2-2)/(P-P\inv)^2>0$, 
whereas $c$ changes its sign at $\mu=0.$
It is easy to verify
$$
\frac89<\mu_a<1.
$$
Furthermore $b$ changes its sign at $\mu=\mu_b:=P/2(P-P\inv)\in(1/2, 2/3).$

\noindent\underline{Case when $\mu\le0$:}
In this case, $a\ge0, c\ge0$ and $b\ge0.$
Since $ 2(1-|c|) - |b|= 2-3P+6P\mu<0,$ Proposition \ref{prop::max} leads to 
$$
\Phi(\mu)
= \frac{1}{3P} \left(a+b+c\right) 
= (1-\mu)P-1.
$$

\noindent\underline{Case when $\mu\ge\mu_a$:}
In this case, $a\le 0, b\le0$ and $c\le0$ and thus $ac\ge0.$
Since $2(1-|c|)-|b|=2+3P-6P\mu<0$ for $\mu\ge\mu_a>1/2+1/3P,$
we have by Proposition \ref{prop::max}
$$
\Phi(\mu)=\frac1{3P}(-a-b-c)=(\mu-1)P+1.
$$

\noindent\underline{Case when $0<\mu<\mu_a$:}
In this case, $a>0, c<0$ and thus $ac<0.$
We compute $b^2+4ac(c^{-2}-1)=H(\mu)/\mu,$ where $H$ is a quadratic
polynomial in $\mu$ given by
$$
H(\mu)=-36\mu^2+(4+P^2+4P^4)\mu-4P^2(P^2-2).
$$
The roots of $H(\mu)$ are given by
$$
\mu_0^\pm=\frac1{72}
\left(4+P^2+4P^4\pm\sqrt{16P^8+8P^6-543P^4+1160P^2+16}\right).
$$
Since $H(2/3)=-2(P^2-4)(2P^2-5)/3<0,$ 
$H(\mu_a)=9P^4(P^2-3)^2(P^2-2)/(P^2-1)^4>0$
and
$H(4/3)=4(P^2-4)(P^2+11)/3>0,$ the roots are real and
satisfy $2/3<\mu_0^-<\mu_a<4/3<\mu_0^+.$
Note that $H(\mu)<0$ for $\mu\in(-\infty,\mu_0^-)\cup(\mu_0^+,+\infty)$
and that $H(\mu)\ge0$ for $\mu\in[\mu_0^-,\mu_0^+].$
Since $2(1-|c|)-|b|=2(1+c)+b=2+(3-6\mu)P<2-P<0$ for $\mu\ge\mu_0^-(>2/3),$
the first case in \eqref{eq::Q} does not occur.

We next analyze the condition $b^2<4(1+|c|)^2,$ which is equivalent
to $|b|<2(1+|c|)=2(1-c)$ in the present case.
We observe that $b<2(1-c)$ precisely when $\mu>\mu_1=1/2-1/3P$
whereas $-b<2(1-c)$ precisely when $\mu<\mu_1':=P(3P+2)/6(P^2-2).$
Note here that $1/3<\mu_1<1/2<\mu_1'<4/3.$
Hence, for $\mu\in(0,\mu_a),$ we see that
$b^2<4(1+|c|)^2$ if and only if $\mu_1<\mu<\mu_1'.$
Hence, by the second case of \eqref{eq::Q}, we obtain
$$
\Phi(\mu)=\frac1{3P}\left(1+a+\frac{b^2}{4(1-c)}\right)
$$
for $\mu_1<\mu<\mu_2=\min\{\mu_0^-,~\mu_1'\}.$
Substituting the explicit forms of $a,b,c,$ we obtain the expression
in the theorem.
Here, keeping $\mu_1'<4/3$ in mind, we see that
$\mu_1'>\mu_0^-$ if and only if
$$
H(\mu_1')=-\frac{P(2P-1)(P^2-4)U(P)}{6(P^2-2)^2}>0,
$$
where $U(P)$ is the quartic polynomial given in Theorem \ref{thm::main}.
One can check that the polynomial $U(P)$ has a unique root $P_*\approx 2.88965$
in the interval $2<P<+\infty.$
Thus $\mu_2=\mu_0^-$ if $2<P\le P_*$ and $\mu_2=\mu_1'$ if $P_*\le P<+\infty.$

When either $0<\mu\le\mu_1$ or $\mu_2\le \mu<\mu_a,$
we have $Y(a,b,c)=R(a,b,c)$ in \eqref{eq::Q}.
We shall take a closer look at these cases.

\medskip
\noindent
{\sl Subcase when $0<\mu<\mu_1$:}
Since $\mu_1<1/2<\mu_b,$ we have $b>0$ in this case.
%Also, noting that
%\begin{align*}
%|b|-4|a|&=b-4a=3P(9-4P^2)+6P\inv(3-5P^2+2P^4)\mu \\
%&< 3P(9-4P^2)+3P\inv(3-5P^2+2P^4)=3P\inv(3+4P^2-2P^4)<0
%\end{align*}
%for $2<P,$ we see that the second case in \eqref{eq::R} does not occur.
%We next examine the condition for the first case in \eqref{eq::R}.
We compute
$$
|ab|-|c|(|b|+4|a|)=ab+c(b+4a)=9\big[2P^2(P^2-1)\mu^2-(3P^4-4P^2-2)\mu+P^2(P^2-2)
\big].
$$
Note that the above quadratic polynomial in $\mu$ is convex and
has the axis of symmetry $\mu=(3P^4-4P^2-2)/4P^2(P^2-1)>1/2>\mu_1.$
Therefore, it is decreasing in $0<\mu<\mu_1$ and thus
\begin{align*}
|ab|-|c|(|b|+4|a|)&\ge
9\big[2P^2(P^2-1)\mu_1^2-(3P^4-4P^2-2)\mu_1+P^2(P^2-2)\big] \\
&=\frac{9}{2P}(6P^4-5P^3-12P^2+14P-12)>0
\end{align*}
for $P>2.$
Hence, by the first case of \eqref{eq::R} in Proposition \ref{prop::max},
we have $\Phi(\mu)=R(a,b,c)/3P=(a+b+c)/3P=(1-\mu)P^2-1.$

\medskip
\noindent
{\sl Subcase when $\mu_2<\mu<\mu_a$:}
First note that $\mu_1'-\mu_b=P(P+2)(2P-1)/6(P^2-1)(P^2-2)>0.$
We also have $\mu_b<2/3<\mu_0^-.$
Thus, we observe that $\mu_b<\mu_2=\min\{\mu_0^-,\mu_1'\},$
which implies that $b<0$ in this case.
Therefore, $|ab|-|c|(|b|+4|a|)=-ab+c(-b+4a)=-9P^{-2}F(\mu),$ where
$$
F(\mu)=2(P^2-1)(P^2-2)^2\mu^2-P^2(3P^4-12P^2+14)\mu+P^4(P^2-2).
$$
The discriminant of $F(\mu)$ is $D=P^4V(P),$ where $V(P)$ is given
in Theorem \ref{thm::main}.
One can see that the polynomial $D$ in $P$
has exactly two roots $P_1, P_2$ in the interval
$2<P<+\infty$ with $P_1\approx 2.05313<P_2\approx 2.82343$
and that $D\ge0$ on $P>2$ if and only if either $2<P\le P_1$ or
$P_2\le P.$
The axis of symmetry of $F(\mu)$ is 
$\mu=\mu_F:=P^2(3P^4-12P^2+14)/4(P^2-1)(P^2-2)^2.$
Since
$$
\mu_F-1
%\frac{P^2(3P^4-12P^2+14)}{4(P^2-1)(P^2-2)^2}-1
=\frac{P^2(-P^6+8P^4-18P^2+16)}{4(P^2-1)(P^2-2)^2}>0\quad (2<P\le 2.2),
$$
we have $F(\mu)>F(1)=2(P^2-4)>0$ for $\mu<1$ and $2<P\le P_1.$
Since $F(\mu)>0$ for all $\mu\in\R$ when $P_1<P<P_2,$ we conclude that
$|ab|-|c|(|b|+4|a|)=-9P^{-2}F(\mu)<0$ for $\mu<\mu_a(<1)$ and $2<P<P_2.$

Solving the equation $F(\mu)=0,$ we write the solutions as
$$
\mu_3^\pm=\frac{P^2(3P^4-12P^2+14)\pm P^2\sqrt{P^8-16P^6+84P^4-176P^2+132}}%
{4(P^2-1)(P^2-2)^2}
$$
for $P\in[P_2,+\infty).$
Note that $F(\mu)>0$ for $\mu\in(-\infty,\mu_3^-)\cup(\mu_3^+,+\infty)$
and that $F(\mu)\le 0$ for $\mu\in[\mu_3^-,\mu_3^+].$
As above, we compute
$$
\mu_a-\mu_F
%\mu_a-\frac{P^2(3P^4-12P^2+14)}{4(P^2-1)(P^2-2)^2}
=\frac{P^2(P^6-9P^4+22P^2-18)}{4(P^2-1)^2(P^2-2)^2}>0\quad (2.5<P),
$$
and
$$
F(\mu_a)=\frac{P^4(P^2-2)(P^2-3)}{(P^2-1)^3}>0,
$$
both of which imply that $\mu_3^+<\mu_a$ for $P_2\le P.$
On the other hand, for $2<P,$ we see that
$$
F(\mu_1')
=-\frac{P^2(P-2)(6P^4-P^3-38P^2-28P+4)}{18(P^2-2)}
=-\frac{P^2(P-2)U(P)}{18(P^2-2)}
\le 0
$$
if and only if $P_*\le P,$ where $P_*$ is the unique root of
$U(P)$ in $2<P<+\infty$ as was introduced above.
Hence, $\mu_3^-\le \mu_1'=\mu_2\le \mu_3^+$ when $P_*\le P,$ and
either $\mu_1'<\mu_3^-$ or $\mu_3^+<\mu_1'$ when $P_2\le P<P_*.$
In view of the fact that 
\begin{align*}
(\mu_3^--\mu_1')\Big|_{P=P_2}
&=\frac{P_2^2(3P_2^4-12P_2^2+14)}{4(P_2^2-1)(P_2^2-2)^2}
-\frac{P_2(3P_2+2)}{6(P_2^2-2)} \\
&=\frac{P_2(3P_2^5-4P_2^4-18P_2^3+12P_2^2+30P_2-8)}{12(P_2^2-1)(P_2^2-2)^2}
\approx 0.049>0,
\end{align*}
we can conclude, by continuity, that $\mu_2=\mu_1'<\mu_3^-$ for $P_2\le P<P_*.$
(In particular, we see that $\mu_0=\mu_1'=\mu_3^-$ at $P=P_*.$
Look around the point $(P_*,\mu_2(P_*))$ in Figure \ref{fig::Pmu}.
We wonder if this is just an incidence.)

Similarly, we have $|c|(|b|-4|a|)-|ab|=-c(-b-4a)+ab=9P\inv G(\mu),$ where
$$
G(\mu)=2P^2(P^2-1)\mu^2-(3P^4-4P^2-2)\mu+P^2(P^2-2).
$$
Solving the equation $G(\mu)=0,$ we write the solutions as
$$
\mu_4^\pm=\frac{3P^4-4P^2-2\pm\sqrt{P^8-12P^4+16P^2+4}}{4P^2(P^2-1)},\quad
2<P.
$$
Here, we note that $P^8-12P^4+16P^2+4=(P^4-6)^2+16P^2-32>132$ for $2<P.$
We now compute $G(\mu_a)=P^2(P^-2)(P^2-3)/(P^2-1)^3>0.$
Since the axis $\mu=\mu_G:=(3P^4-4P^2-2)/4P^2(P^2-1)$ of $G(\mu)$ satisfies
$\mu_G<3/4<\mu_a,$ we have $\mu_4^+<\mu_a.$
On the other hand, since
$$
\mu_4^--\frac12=
\frac{-P^4-2P^2-2-\sqrt{P^8-12P^4+16P^2+4}}{4P^2(P^2-1)}<0,
$$
we get $\mu_4^-<1/2<\mu_2$ for $2<P.$
We now show that $\mu_0^-<\mu_4^+$ for $2<P,$ from which the inequality
$\mu_2<\mu_4^+$ will follow.
Since
$16P^8+8P^6-543P^4+1160P^2+16-(4P^4-8P^2-8)^2=3(P^2-4)(24P^4-85P^2+4)>0,$
we have
$$
72\mu_0^+>4+P^2+4P^4+\sqrt{(4P^4-8P^2-8)^2}
=8P^4-7P^2-4>0
$$
for $P>2.$
Therefore,
$$
\mu_0^-=\frac{P^2(P^2-2)}{9\mu_0^+}<\frac{8P^2(P^2-2)}{8P^4-7P^2-4}.
$$
On the other hand, since $P^8-12P^4+16P^2+4=(P^4-6)^2+16(P^2-2)>(P^4-6)^2,$
we obtain
$$
\mu_4^+>\frac{3P^4-4P^2-2+(P^4-6)}{4P^2(P^2-1)}
=\frac{(P^2+1)(P^2-2)}{P^2(P^2-1)}.
$$
Because
$$
\frac{(P^2+1)(P^2-2)}{P^2(P^2-1)}-\frac{8P^2(P^2-2)}{8P^4-7P^2-4}
=\frac{(P^2-2)(9P^4-11P^2-4)}{P^2(P^2-1)(8P^4-7P^2-4)}>0
$$
for $P>2,$ the inequality $\mu_0^-<\mu_4^+$ follows as required.

We now summarize the above observations as follows.
Let $D=\{(P,\mu): 2<P, \mu_2(P)<\mu<\mu_a(P)\}.$
Here, we write $\mu_2$ {\it etc.}~as functions of $P.$
We divide $D$ into three parts $D_1, D_2, D_3$ according as
the first, second, third case occurs in \eqref{eq::R}, respectively.
Then,
$D_1=\{(P,\mu): P_2\le P<P_*, \mu_3^-(P)\le \mu\le\mu_3^+(P)\}
\cup\{(P,\mu): P_*\le P, \mu_2(P)<\mu\le\mu_3^+(P)\}$ and
$D_2=\{(P,\mu): \mu_4^+(P)\le \mu<\mu_a(P)\}.$
Since $D_1$ and $D_2$ are disjoint, we have necessarily that
$\mu_3^+<\mu_4^+$ for $P_2\le P.$
Note here that $\Phi(\mu)=(a-b+c)/3P=P^2-3-\mu(P^2-2)^2/P^2$
for $(P,\mu)\in D_1,$
that  $\Phi(\mu)=(-a-b-c)/3P=(\mu-1)P^2+1$
for $(P,\mu)\in D_2$ and that
\begin{align*}
\Phi(\mu)&=\frac{(a-c)}{3P}\sqrt{1-\frac{b^2}{4ac}} \\
&=(1-\mu)P(P^2-2)\sqrt{\frac{P^2-4\mu}{4\mu\{(1-\mu)(P^2-1)^2-1\}}}
\end{align*}
for $(P,\mu)\in D_3.$

Finally, setting $\mu_4=\mu_4^+$ for simplicity, we complete
the proof of Theorem \ref{thm::main}.

\section{Variability region of $a_3-a_2^2$}

We first note that the class $\Co_p$ is not rotationally invariant
for $0<p<1$ due to the presence of a pole at $p.$
It is therefore more natural to consider the variability region of
the Fekete-Szeg\H{o} functional $\Lambda_\mu$ over $\Co_p$
rather than its modulus only.
The present section will be devoted to the study of the variability
region of $\Lambda_1(f)=a_3-a_2^2$ because of its importance.
Let
$$
W_p=\{\Lambda_1(f): f\in\Co_p\}
$$
for $0<p<1.$

\begin{figure}
\begin{center}
\includegraphics[width=0.5\textwidth]{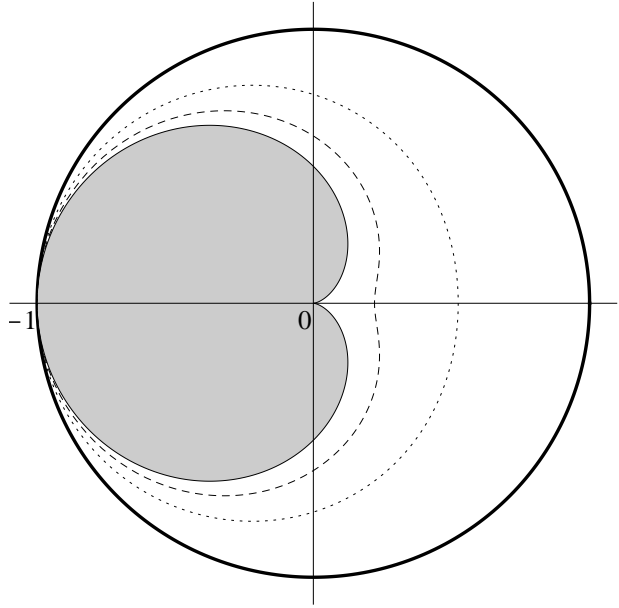}
\caption{A couple of $\Omega_p$'s (the inside of dotted and dashed curves), 
the intersection cardioid and the unit disk}
\label{fig:card}
\end{center}
\end{figure}

In the following, we fix $p\in(0,1)$ and put $P=p+1/p>2.$
Let
$$
f_\zeta(z)= \frac{z-T_p(p \zeta)z^2}{(1-z/p)(1+pz)}
= \sum_{n=1}^{\infty} \frac{1-p^{2n}\zeta}{p^{n-1}(1-p^2\zeta)} z^n  
= \sum_{n=1}^{\infty} A_n(\zeta) z^n
$$
for $z\in\D$ and $\zeta\in\bD.$
Here, $T_p$ is defined in \eqref{eq::Ta}.
One can check that $f_\zeta$ belongs to $\Co_p$ and corresponds to 
$\varphi(z)=T_p(\zeta T_p(z))$ through \eqref{equ::intp}.
As Avkhadiev and Wirths \cite{AW07} pointed out, the function $f_\zeta$
with $|\zeta|=1$ is extremal in important problems for the class $\Co_p.$
Indeed, they proved that the closed disk $A_n(\bD)$ is the variability region
of the coefficient functional $a_n(f)$ for $f\in\Co_p$
(see also \cite{OS15}).
We now compute
$$
\Lambda_1(f_\zeta)=A_3(\zeta)-A_2(\zeta)^2
= -\frac{(1-p^2)^2\zeta}{(1-p^2\zeta)^2} 
%= - \frac{(1-p^2)^2}{p^2} K(p^2 \zeta),
= - (P^2-4) K(p^2 \zeta),
$$
where $K(z)=z/(1-z)^2$ is the Koebe function.
One might expect that the image
$$
\Omega_p=\{-(P^2-4)K(p^2 \bar z): |z|\le 1\}
$$
would coincide the variability region $W_p.$
By accident, the form of $A_3-A_2^2$ is same as the second Hankel determinant
$a_2a_4-a_3^2$ of order $2$ for $f_\zeta;$ namely,
$$
A_2(\zeta)A_4(\zeta)-A_3(\zeta)^2
= -\frac{(1-p^2)^2\zeta}{(1-p^2\zeta)^2} 
= A_3(\zeta)-A_2(\zeta)^2.
$$
The authors investigated in \cite{OS15}
the set $\Omega_p$ in the context of the second
Hankel determinant and found that 
$\Omega_p\subset\Omega_q$ for $0<q<p<1$ and that
$$
\bigcup_{0<p<1}\Omega_p=\D\cup\{-1\}
\aand
\bigcap_{0<p<1}\Omega_p=\{-(1+z)^2/4: |z|\le 1\}.
$$
Note that $\{-(1+z)^2/4: |z|\le 1\}$ is a closed Jordan domain
bounded by a cardioid (see Figure \ref{fig:card}).
We also observed in \cite{OS15} that the variability region of $a_2a_4-a_3^2$
for $\Co_p$ is properly larger than $\Omega_p.$
In the case of $a_3-a_2^2,$ rather surprisingly, the expected result
partially holds and a phase transition occurs.

\begin{Theorem}
Let $0<p<1.$
The variability region $W_p$ of $a_3-a_2^2$ for $\Co_p$ satisfies
$\Omega_p\subset W_p\subset\bD.$
Moreover, $W_p=\Omega_p$ for $0<p\le p_0$ and $W_p\ne\Omega_p$
for $p_0<p<1,$ where
$$
p_0=\frac{1+\sqrt{37}-\sqrt{2(1+\sqrt{37})}}6
\approx 0.553175.
$$
\end{Theorem}

\begin{pf}
Letting $\sigma=T_{p^2}(\zeta),$ we have the representation
$$
A_3(\zeta)-A_2(\zeta)^2
=-\frac{p^2}{(1+p^2)^2}\left[1-(p^2+p^{-2})\sigma+\sigma^2\right]
=-P^{-2} h(\sigma),
$$
where
$$
h(\sigma)=1-t\sigma+\sigma^2,\quad t=P^2-2>2.
$$
Hence, $\Omega_p=-P^{-2}h(\bD).$
One can easily check that $h$ is univalent on $\bD.$
Let $\Delta_r$ be the image of the closed disk $|z|\le r$ under the
mapping $h$ for $0\le r\le 1.$
For $\zeta,\omega\in\partial\D,$ the sharp inequality
$$
|h(\zeta)-h(r\omega)|=|\zeta-r\omega||\zeta+r\omega-t|
\ge (1-r)(t-1-r)=h(r)-h(1),
$$
holds.
Hence, the Euclidean distance $\delta_r$ between $\partial\Delta_r$ and
$\partial\Delta_1$ is given as $(1-r)(t-1-r)=(1-r)(P^2-3-r)$
for $0\le r\le1.$
Note that if $|w-h(\sigma)|\le\delta_r$ for some $\sigma\in\C$
with $|\sigma|=r,$ then $w\in\Delta_1.$

Letting $\mu=1$ in \eqref{eq::FS}, we obtain
the following representation of $\Lambda_1(f)$ for $f(z)=z+a_2z^2+a_3z^3
+\cdots\in\Co_p:$
\begin{align*}
a_3-a_2^2&=-P^{-2}+(2P^{-2}-1)\sigma_0-P^{-2}\sigma_0^2
-\frac{(1-|\sigma_0|^2)\sigma_1}{3P} \\
&=-P^{-2}\left[
h(-\sigma_0)+(1-|\sigma_0|^2)\sigma_1 P/3
\right]
\end{align*}
for some $\sigma_0,\sigma_1\in\bD.$
Put $r=|\sigma_0|.$
Then $h(-\sigma_0)\in\partial\Delta_r.$
If $(1-r^2)P/3\le \delta_r,$ we have $a_3-a_2^2\in-P^{-2}\Delta_1=\Omega_p.$
Since
$$
\delta_r-(1-r^2)P/3
=\frac{1-r}{3}\left[3P^2-(1+r)P-9-3r\right]
\ge\frac{1-r}{3}\left[3P^2-2P-12\right],
$$
we have $(1-r^2)P/3\le\delta_r$ for $P\ge P_0:=(1+\sqrt{37})/3
\approx 2.36092,$
which is the larger root of the polynomial $3P^2-2P-12.$
Note that $p_0$ is determined by $P_0=p_0+1/p_0.$
Thus we have shown that $W_p\subset\Omega_p$ for $0<p\le p_0.$

We next assume that $2<P<P_0.$
Since $3P^2-2P-12<0,$ we can find an $r\in(0,1)$ such that
$h(r)-h(1)-(1-r^2)P/3=(1-r)[3P^2-(1+r)P-9-3r]/3<0.$
We choose $\sigma_0=-r$ and $\sigma_1=1.$
Then there is a function $f(z)=z+a_2z^2+a_3z^3+\cdots$ in $\Co_p$
satisfying \eqref{eq::FS} with $\mu=1:$
$$
a_3-a_2^2
=-P^{-2}\left[h(r)+(1-r^2)P/3\right].
$$
Therefore, we get
$$
a_3-a_2^2=-P^{-2}\left[h(1)+\{h(r)-h(1)+(1-r^2)P/3\}\right]
>-P^{-2}h(1)=1-4P^{-2},
$$
which implies that $a_3-a_2^2\in W_p\setminus\Omega_p$
because $\Omega_p\cap\R
%=[-P^{-2}h(-1),-P^{-2}h(1)]
=[-1,1-4P^{-2}].$
\end{pf}

%\bibliography{papers}

\def\cprime{$'$} \def\cprime{$'$} \def\cprime{$'$}
\providecommand{\bysame}{\leavevmode\hbox to3em{\hrulefill}\thinspace}
\providecommand{\MR}{\relax\ifhmode\unskip\space\fi MR }
% \MRhref is called by the amsart/book/proc definition of \MR.
\providecommand{\MRhref}[2]{%
  \href{http://www.ams.org/mathscinet-getitem?mr=#1}{#2}
}
\providecommand{\href}[2]{#2}

\end{document}